\documentclass[letterpaper, 10pt, journal,twoside]{IEEEtran}  

\IEEEoverridecommandlockouts                              

\usepackage{amsmath,amssymb,color,cite}

\usepackage{hyperref}
\usepackage{xcolor}
\hypersetup{
	colorlinks,
	linkcolor={red!80!black},
	citecolor={blue!80!black},
	urlcolor={blue!80!black}
}

\DeclareMathOperator{\im}{im}

\DeclareMathOperator{\inte}{int}
\DeclareMathOperator{\cl}{cl}

\DeclareMathOperator{\dom}{dom}
\DeclareMathOperator{\graph}{gr}
\DeclareMathOperator{\gr}{gr}
\DeclareMathOperator{\lin}{lin}
\DeclareMathOperator{\Lin}{Lin}

\let\leq\leqslant
\let\geq\geqslant
\let\emptyset\varnothing

\newcommand{\calC}{\ensuremath{\mathcal{C}}}

\newcommand{\calF}{\ensuremath{\mathcal{F}}}

\newcommand{\calK}{\ensuremath{\mathcal{K}}}

\newcommand{\calN}{\ensuremath{\mathcal{N}}}

\newcommand{\calR}{\ensuremath{\mathcal{R}}}
\newcommand{\calS}{\ensuremath{\mathcal{S}}}
\newcommand{\calT}{\ensuremath{\mathcal{T}}}

\newcommand{\calW}{\ensuremath{\mathcal{W}}}

\newcommand{\calY}{\ensuremath{\mathcal{Y}}}

\newcommand{\barr}{\ensuremath{\bar{r}}}

\newcommand{\barx}{\ensuremath{\bar{x}}}

\newcommand{\barH}{\ensuremath{\bar{H}}}

\newcommand{\bmat}{\begin{matrix}}
\newcommand{\emat}{\end{matrix}}
\newcommand{\bbm}{\begin{bmatrix}}
\newcommand{\ebm}{\end{bmatrix}}
\newcommand{\bbma}{\begin{bmatrix*}[r]}
\newcommand{\ebma}{\end{bmatrix*}}
\newcommand{\bpm}{\begin{pmatrix}}
\newcommand{\epm}{\end{pmatrix}}
\newcommand{\bpma}{\begin{pmatrix*}[r]}
\newcommand{\epma}{\end{pmatrix*}}
\newcommand{\bvm}{\begin{vmatrix}}
\newcommand{\evm}{\end{vmatrix}}
\newcommand{\bse}{\begin{subequations}}
\newcommand{\ese}{\end{subequations}}
\newcommand{\beq}{\begin{equation}}
\newcommand{\eeq}{\end{equation}}

\newcommand{\ben}{\begin{enumerate}}
\newcommand{\een}{\end{enumerate}}
\newcommand{\beni}{\begin{enumerate}}
\newcommand{\eeni}{\end{enumerate}}
\newcommand{\bena}{\begin{enumerate}}
\newcommand{\eena}{\end{enumerate}}

\newcommand{\bit}{\begin{itemize}}
\newcommand{\eit}{\end{itemize}}
\newcommand{\bthe}{\begin{theorem}}
\newcommand{\ethe}{\end{theorem}}
\newcommand{\blem}{\begin{lemma}}
\newcommand{\elem}{\end{lemma}}
\newcommand{\bprop}{\begin{proposition}}
\newcommand{\eprop}{\end{proposition}}
\newcommand{\bex}{\begin{example}}
\newcommand{\eex}{\end{example}}
\newcommand{\bas}{\begin{assumption}}
\newcommand{\eas}{\end{assumption}}
\newcommand{\bre}{\begin{remark}}
\newcommand{\ere}{\end{remark}}
\newcommand{\bcor}{\begin{corollary}}
\newcommand{\ecor}{\end{corollary}}
\newcommand{\bdfn}{\begin{definition}}
\newcommand{\edfn}{\end{definition}}
\newcommand{\bcon}{\begin{conjecture}}
\newcommand{\econ}{\end{conjecture}}

\newcommand{\inv}{\ensuremath{^{-1}}}
\newcommand{\nonempty}{\ensuremath{\neq\emptyset}}
\newcommand{\pset}[1]{\ensuremath{\{#1\}}}

\newcommand{\zset}{\ensuremath{\pset{0}}}

\newcommand{\R}{\ensuremath{\mathbb R}}

\newcommand{\N}{\ensuremath{\mathbb N}}

\newcommand{\qand}{\quad\text{ and }\quad}

\newcommand{\rrm}{\calR_-}
\newcommand{\rrp}{\calR_+}

\newtheorem{theorem}{Theorem}
\newtheorem{example}{Example}
\newtheorem{remark}{Remark}
\newtheorem{proposition}{Proposition} 
\newtheorem{lemma}{Lemma}
\newtheorem{corollary}{Corollary}
\newtheorem{dfn}{Definition}
\newtheorem{assumption}{Assumption}

\title{On reachability and null-controllability of nonstrict convex processes
}

\author{J. Eising$^{1}$ and M.K. Camlibel$^{2}$
\thanks{$^{1}$ Bernoulli Institute for Mathematics, Computer Science, and Artificial Intelligence, University of Groningen, Nijenborgh 9, 9747 AG, Groningen, The Netherlands, {\footnotesize{\tt j.eising@rug.nl}}}%
\thanks{$^{2}$ Bernoulli Institute for Mathematics, Computer Science, and Artificial Intelligence, University of Groningen, Nijenborgh 9, 9747 AG, Groningen, The Netherlands, {\footnotesize{\tt m.k.camlibel@rug.nl}}}
}

\newcounter{todocounter}
\usepackage[prependcaption,colorinlistoftodos]{todonotes}

\begin{document}

\begingroup

\maketitle
\thispagestyle{empty}
\pagestyle{empty}

\begin{abstract}This paper studies reachability and null-controllability for difference inclusions involving convex processes. Such difference inclusions arise, for instance, in the study of linear discrete-time systems whose inputs and/or states are constrained to lie within a convex cone. After developing a geometric framework for convex processes relying on invariance properties, we provide necessary and sufficient conditions for both reachability and null-controllability in terms of the spectrum of dual processes.
\end{abstract}



\section{Introduction}

\IEEEPARstart{T}{he} motivation behind this paper stems from the question when a constrained linear discrete-time system of the form
\bse\label{eq:system +constraints}
\begin{align}
x_{k+1}=& Ax_k+Bu_k \\
&Cx_k+Du_k\in\calY
\end{align}
\ese
is reachable or null-controllable. Although this is a natural and fundamental question, it is still open in full generality. Indeed, all existing results in the literature deal with particular cases, which are summarized in Table~\ref{tbl:literature}. The references cited in Table~\ref{tbl:literature} all provide necessary and sufficient spectral conditions for the problems they study. Similar results were also provided in continuous-time for which we refer the reader to \cite{Br:72,So:84,Ng:85,HC:07}.

In this paper, we focus on convex conic constraints, i.e. $\calY$ is a convex cone. This means that we will work in the more general framework of difference inclusions of the form
\begin{equation}\label{eq:diffinc intro}
x_{k+1}\in H(x_k)
\end{equation}
where $H$ is a convex process, a set-valued map whose graph is a convex cone (see Section~\ref{sec:Conv proc}). The constrained system \eqref{eq:system +constraints} can be represented in the form \eqref{eq:diffinc intro} by taking $H(x) = \{Ax+Bu \mid Cx+Du\in \calY\}$.

Systems of the form \eqref{eq:diffinc intro} are encountered in various contexts. Examples include von Neumann-Gale economic growth models \cite{makarov:77}, cable-suspended robots \cite{j37,oh:05-2} and chemical reaction networks \cite{Angeli:09}. 

One of the advantages of this framework is the ease of studying reachability and null-controllability problems by employing invariance properties. The notions of reachability and null-controllability are not only of interest from a theoretical 
\begin{table}[h]
	\begin{center}
		\begin{tabular}{|r|c|c|p{5.5cm}|}
			\hline
			Ref. & R & N & extra conditions\\
			\hline
			\hline
			\cite{Trentelman:01} & \checkmark & \checkmark & $\calY$ is a subspace\\
			\hline
			\hline
			\cite{So:84} & & \checkmark & $C=0$, $D=I$, $\calY$ is a bounded convex set containing the origin in its interior\\
			\hline
			\hline
			\cite{Ev:85}& \checkmark & & $C=0$, $D=I$, $\calY$ is a closed convex cone\\
			\hline
			\hline
			\cite{Ng:86}& \checkmark & \checkmark & $C=0$, $D=I$, $\calY$ is a hyperbolic convex set containing the origin \\
			\hline
			\hline
			\cite{HC:08} & & \checkmark & $D+C(sI-A)\inv B$ is right invertible and $\calY$ is a solid polyhedral cone\\
			\hline
			\hline
			\cite{HC:08} & & \checkmark & $D+C(sI-A)\inv B$ is right invertible as a rational matrix and $\calY$ is a bounded convex set containing the origin in its interior \\
			\hline
			\hline
			\cite{j37} & \checkmark & & $\calY$ is a convex cone such that $\im\bbm C & D\ebm \cap \calY$ is solid and $(\im D+C\calT^*)\cap\inte(\calY)\neq\emptyset$\\
			\hline
			\hline
			\cite{j37} & \checkmark & & $\calY$ is a closed pointed convex cone such that $\im\bbm C & D\ebm \cap \calY$ is solid and $(\im D+C\calT^*)\cap\calY=\zset$\\
			\hline
			\hline
			\cite{Kaba-thesis} & \checkmark & & $\calY$ is a hyperbolic convex set such that $\im\bbm C & D\ebm \cap \calY$ is solid and $\im D+C\calT^*+\calY=\calY-\calY$\\
			\hline
		\end{tabular}
	\end{center}
	\label{tbl:literature}
	\caption{Results on reachability (R) and null-controllability (N)}
	\vspace*{-8mm}
\end{table}
\noindent point of view, but are also fundamental prerequisites for virtually any design problem. For instance, null-controllability plays a role in constrained model predictive control schemes, guaranteeing feasibility of optimal control in the presence of endpoint constraints (see e.g. \cite{Mayne:00}).

To the best of our knowledge, the earliest work studying reachability and null-controllability of difference inclusions is \cite{PD:94}. This paper provides necessary and sufficient conditions for reachability when the convex process is {\em strict\/}, i.e. $H(0)\neq\emptyset$ for all $x$. In the context of constrained systems, strictness corresponds to the restricted case of constraints involving {\em only\/} inputs {\em but not\/} states (see Section~\ref{sec:Conv proc} for a formal definition). As such, it is rather restrictive. The paper \cite{PD:94} also provides necessary and sufficient conditions for null-controllability under the assumption that both the underlying convex process and its inverse are {\em strict\/}. For nonstrict processes, the null-controllability problem has not been studied in the literature yet whereas the only work dealing with the reachability problem is \cite{j37}. For the work on reachability of convex processes in the continuous-time setup, we refer the reader to the seminal paper of \cite{AFO:86} for the strict case and the paper \cite{Seeger:01} studying duality relations for the nonstrict case. 

The contributions this paper brings about are three-fold. Firstly, we develop a novel geometric framework, based on invariance properties, for convex processes. This framework not only paves the road to the study of reachability and null-controllability problems for nonstrict convex processes but also opens up the possibility of extending the geometric approach to linear systems (see e.g. \cite{Wonham1985,Basile1992,Trentelman:01}) towards convex processes. Secondly, we provide necessary and sufficient spectral conditions for reachability under weaker assumptions than those employed in the literature. Thirdly, we provide necessary and sufficient conditions for null-controllability problem for nonstrict convex processes. To the best of our knowledge, this problem has not been studied in the literature before.  

The outline of the paper is as follows. In Section~\ref{sec:prob for}, we formulate the problems of reachability and null-controllability for general difference inclusions. In Section~\ref{sec:Conv proc}, we quickly introduce convex process and related notions. Based on these notions, we are able to provide the main results in Section~\ref{sec:main res}. The proofs of our results require the introduction of several notions and auxiliary results that are presented in Section~\ref{sec:towards}. Section~\ref{sec:proofs} gives the proofs of the main results. Finally, the paper closes with conclusions in Section~\ref{sec:conc}.

\section{Problem formulation}\label{sec:prob for}
A \textit{set-valued map} $H:\R^n\rightrightarrows \R^n$ is a map taking elements of $\R^n$ to subsets of $\R^n$. Consider a system described by a \textit{difference inclusion} of the form:
\begin{equation}\label{eq:diffinc} x_{k+1} \in H(x_k) \quad k\geq 0
 \end{equation}
where $H:\R^n\rightrightarrows \R^n$ is a set-valued map.	
By a \textit{trajectory} of \eqref{eq:diffinc}, we mean a sequence $(x_k)_{k\in\N }$ such that \eqref{eq:diffinc} holds for all $k\geq 0$.

Next, we define a number of sets associated with the system \eqref{eq:diffinc}. The \textit{behavior} (e.g. \cite{Willems:91}) is the set of all trajectories: 
\[ \mathfrak{B}(H) = \left\lbrace (x_k) \in (\mathbb{R}^n)^\mathbb{N}\mid (x_k) \textrm{ is a trajectory of } \eqref{eq:diffinc}\right\rbrace . \]
The \textit{feasible set} is the set of states from which a trajectory emanates:	
	\[ \mathcal{F}(H)= \left\lbrace\xi \mid \exists (x_k)\in \mathfrak{B}(H) \textrm{ with } x_0=\xi \right\rbrace.\] 
We also define the set of $q$-step trajectories as 
\[ \mathfrak{B}_q(H) = \left\lbrace (x_k)_{k=0}^q \in (\mathbb{R}^n)^{q+1}\mid (x_k) \textrm{ satisfies } (\ref{eq:diffinc})\right\rbrace . \] 

In this paper, we are interested in reachability and null-controllability of system \eqref{eq:diffinc}. We define the set of all states that can be reached in finite steps from the origin as the \textit{reachable set} and the set of all states that can be steered in finite steps to the origin as the \textit{null-controllable set}. These sets are, respectively, denoted by $\mathcal{R}(H)$ and $\mathcal{N}(H)$:
\bse\label{eq:reach null def}
\begin{align}
\!\!\calR(H) &\!=\! \big\lbrace \xi \mid \exists (x_k)_{k=0}^q \in\mathfrak{B}_q(H) \textrm{ s.t. } x_0=0, x_q=\xi \big\rbrace \label{eq:reach def}\\ 
\!\!\!\calN(H) &\!=\! \big\lbrace \xi \mid \exists (x_k)_{k=0}^q \in\mathfrak{B}_q(H) \textrm{ s.t. } x_0=\xi, x_q=0 \big\rbrace\label{eq:null def}
\end{align}
\ese
We say the system \eqref{eq:diffinc} is reachable (null-controllable) if every feasible state is reachable (null-controllable), that is, $\calF(H)\subseteq\calR(H)$ $\big(\calF(H)\subseteq\calN(H)\big)$. 	

%
%
%

In this paper, we will focus on a specific class of set-valued maps, namely the so-called convex processes, and derive necessary and sufficient conditions for both reachability and null-controllability. 

%
%
	
\section{Convex processes}\label{sec:Conv proc}
A convex cone is a nonempty convex set that is closed under nonnegative scalar multiplication. For two nonempty convex sets $\calS,\calT\subseteq\R^n$ and scalar $\rho\in\R$, we define the (Minkowski) sum and scalar product of sets as: 
	\[ \calS+\calT = \{s+t \mid s\in\calS, t\in\calT\}, \quad \rho \calS = \{\rho s \mid s\in\calS\}. \] 

%

A set-valued map $H:\R^n\rightrightarrows \R^n$ is called a {\em convex process\/}, a {\em linear process\/}, {\em closed\/} if its graph
\[ \graph(H) = \{(x,y)\in \R^n\times \R^n \mid y\in H(x)\}\]
is a convex cone, a subspace, closed, respectively.
	 
\begin{example}
Consider the following linear system with constraints:
\begin{equation}\label{eq:dynamics} x_{k+1} = Ax_k+Bu_k, \quad y_k= Cx_k+Du_k\in \mathcal{Y}\end{equation}
where $A\in\R^{n\times n}$, $B\in\R^{n\times m}$, $C\in\R^{p\times n}$, $D\in\R^{p\times m}$, and $\calY\subseteq\R^p$. Consider in addition the set-valued map
\begin{equation}\label{eq:lin const proc}
 H(x) = \{ Ax+Bu\mid Cx+Du\in \mathcal{Y} \}.
\end{equation}
It can be readily verified that $H$ is a convex (linear) process  if $\calY$ is a convex cone
(subspace).	Clearly, there is one-to-one correspondence between the state trajectories of   \eqref{eq:dynamics} and those of \eqref{eq:diffinc}. 
\end{example}

The \textit{domain} and \textit{image} of $H$ are defined as $\dom(H)= \{x\in \R^n\mid H(x)\neq \emptyset \}$ and $\im (H) = \{y\in \R^n : \exists\, x \textrm{ s.t. } y\in H(x)\}$. If $\dom(H)=\R^n$, we say $H$ is \textit{strict}. The {\em inverse\/} of $H$ is defined by:
$$
(y,x)\in\gr(H\inv)\iff (x,y)\in\gr(H).
$$
Clearly, one has $\dom(H\inv)=\im(H)$ and vice versa. 
	
For a convex cone $ \calK\subseteq\R^n$, we define $\lin(\calK) =-\calK\cap \calK$ and $\Lin(\calK)=\calK-\calK$. Note that $\lin(\calK)$ is the largest subspace contained in $\calK$ whereas $\Lin(\calK)$ is the smallest subspace that contains $\calK$. 
	
Let $H$ be a convex process. Associated with $H$, we define two linear processes $L_-$ and $L_+$ by
	\beq\label{eq:def of L- and L+}
	\graph(L_-)=\lin\big(\gr(H)\big)\textrm{ and } \graph(L_+)=\Lin\big(\graph(H)\big).
	\eeq
By definition, we therefore have
	\beq\label{eq:lh-+}
\graph(L_-)\subseteq\graph(H)\subseteq\graph(L_+).
	\eeq
It is clear that $L_-$ and $L_+$ are, respectively, the largest and the smallest (with respect to the graph inclusion) linear processes satisfying \eqref{eq:lh-+}. We call $L_-$ and $L_+$, respectively, the minimal and maximal linear processes associated with $H$. If $H$ is not clear from context, we write $L_-(H)$ and $L_+(H)$ in order to avoid confusion.	

\begin{example}
	Let $H$ be of the form \eqref{eq:lin const proc}. It can be shown that if the set $\{u\mid Bu=0,Du\in\calY\}$ is a subspace, then
	\[ L_- = \{ Ax+Bu \mid Cx+Du\in \lin(\mathcal{Y}) \}. \] 
	Furthermore, if $\im \begin{bmatrix} C & D \end{bmatrix} +\calY $ is a subspace, then
	\[ L_+ = \{ Ax+Bu \mid Cx+Du\in \Lin(\mathcal{Y}) \}. \] 
\end{example}
		
For a nonempty set $\calC\subseteq \R^n $, we define the \textit{negative} and \textit{positive  polar cone}, respectively,
\begin{align*}
\calC^- &= \{ y\in\R^n \mid \langle x,y\rangle\leq 0 \hspace{1em} \forall x\in \calC\}, \\
\calC^+ &= \{ y\in\R^n \mid \langle x,y\rangle\geq 0 \hspace{1em} \forall x\in \calC\}.
\end{align*}

Based on these, we define \textit{dual} processes $H^-$ and $H^+$ of $H$ as follows:
\bse
\begin{align} \label{eq:def of dual} 
  p\in H^-(q) &\iff \langle p,x\rangle\geq \langle q,y\rangle\quad \forall\, (x,y)\in\gr(H)  \\
		p\in H^+(q) & \iff \langle p,x\rangle\leq \langle q,y\rangle\quad \forall\, (x,y)\in\gr(H) 
\end{align}
\ese 
\begin{example}
	Let $H$ again be of the form $\eqref{eq:lin const proc}$. It can be shown that if $\im \begin{bmatrix} C & D \end{bmatrix} +\calY $ is a subspace, then
	\[ H^-(x)= \left\{ A^\top x + C^\top u \mid u\in \calY^+, B^\top x+D^\top u=0\right\}  \] 
\end{example} 

We close this section with the definition of eigenvalues. We say $\lambda\in \R$ is an {\em eigenvalue\/} of $H$ if there exists a nonzero vector $\xi\in\R^n $ such that $\lambda \xi \in H(\xi)$. Such a vector $\xi$ will be called an eigenvector corresponding to $\lambda$.

\section{Main results}\label{sec:main res}

In this section, we provide the main results of the paper whose proofs can be found in Section~\ref{sec:proofs}. We begin with reachability. For the sake of brevity, we define $\rrm = \calR(L_-)$ and $\rrp= \calR(L_+)$ in the sequel. 
		
\begin{theorem}\label{thm:reach}
	Let $H$ be a convex process such that $\dom H +\rrm =\R^n$. Then, $H$ is reachable if and only if $\calR_+=\R^n$ and $H^{-}$ has no nonnegative eigenvalues. 	
\end{theorem}
The advantage of this result over the existing results of \cite[Thm. 6.3]{j37} will be illustrated by an example.
\begin{example} 
Consider a constrained linear system of the form \eqref{eq:dynamics} where
\[A=1,\, B=1,\, C=\bbm1\\0\ebm,\, D=\bbm0\\1\ebm,\,\text{ and }\calY=\R_+\times \R .
\]
One can equivalently consider the difference inclusion \eqref{eq:diffinc} where $H$ is of the form \eqref{eq:lin const proc}. Note that $L_-(x)= \{ x+u \mid x\in \zset, u\in\R \}$. As such, it is immediate that $\rrm = \R$ and thus $\dom(H)+\rrm = \R$. Applying Theorem~\ref{thm:reach}, one can verify that $H$ is reachable. However, \cite[Thm.
6.3]{j37} cannot be employed to infer reachability for this example. Indeed, the assumption $\im D+C\calT^*+\calY=\R^2$ where $\calT^*$ is 
the smallest conditioned invariant subspace (see e.g. \cite{Trentelman:01}) is not satisfied since $\calT^*=\zset$. 
\end{example}

Next, we will study null-controllability. Clearly, we have $\calN(H)=\calR(H\inv)$ from \eqref{eq:reach null def}. Based on this observation, one could try to apply Theorem~\ref{thm:reach} to the convex process $H\inv$. However, reachability of $H\inv$ means that $\calF(H^{-1})\subseteq\calR(H^{-1})=\calN(H)$ whereas null-controllability of $H$ amounts to $\calF(H)\subseteq\calN(H)$. Therefore, Theorem~\ref{thm:reach} cannot be employed to test for null-controllability. For the sake of brevity, we define $\calN_-=\calN(L_-)$ in the sequel. 

\begin{theorem}\label{thm:nullc}
Let $H$ be a convex process such that $\dom H +\rrm = \rrp = \im H+\calN_-=\R^n$. Then, $H$ is null-controllable if and only if $H^{-}$ has no positive eigenvalues. 
\end{theorem}	

The following example illustrates the role played by the assumption $\im H+\calN_-=\R^n$. 
\begin{example}
Let $H:\R\rightrightarrows\R $ be the convex process defined by:
\[ H(x) = \begin{cases} [0,\infty)  &\text{ if } x=0,\\ (0,\infty) &\text{ if } x\neq 0.\end{cases}
\] 

Clearly, $H$ is strict and $\rrm = \{0\}$ and $\rrp=\R$. However, $H$ is not null controllable since $\calN(H) = \{0\}$. 

Let $\barH$ be the closure of $H$, that is, $\bar{H}(x) = [0,\infty)$ for any $x\in \R$. Then the process $\barH$ is also strict. Note that $\calR(L_-(\barH)) = \{0\}$ and $\calR(L_+(\barH))=\R$. Since $0\in\barH(x)$ for every $x\in\R$, the convex process $\barH$ is null-controllable. 

Even though $H$ and $\bar{H}$ have the same dual $H^-$, the former is not null-controllable whereas the latter is. This reveals the role played by the assumption on the image. Indeed, we have $\im H+\calN_-=[0,\infty)$ whereas $\im\barH+\calN(L_-(\barH))=\R$. 
\end{example}
An added benefit of this example is that $\bar{H}$ also is a null-controllable process, where $\im\bar{H}\neq \R$. Therefore it is not covered by the result of \cite{Ng:86}.

\section{Towards the Proofs}\label{sec:towards}
Before proving the main results, we need a bit of preparation that will be presented in this section.
\subsection{Convex cones}
Let $\calC$ be a convex cone. Then, the negative polar $\calC^-$ is always closed and moreover $\calC^-=(\cl(\calC))^-$ and $(\calC^-)^-=\cl(\calC)$ where $\cl$ denotes the closure of a set. 
For the sums and intersections of two convex cones $\calC$ and $\calS$, it holds that:
	\begin{equation}\label{eq:dual of sum and intersection} (\calC+\calS)^-  =\calC^-\cap \calS^-, \quad (\calC\cap\calS)^- = \cl (\calC^-+\calS^-). \end{equation}
All the properties of negative polar we mentioned above also hold for the positive polar. 

Note that the sum of two closed convex sets is not necessarily closed. As such, the $\cl$ operator cannot be dropped from the second identity in \eqref{eq:dual of sum and intersection} in general. The following result provides a sufficient condition for the sum of two closed convex cones to be closed as well. 
	\begin{lemma}\label{lemm:sum of closed sets} Let $\calK_1,\calK_2$ be closed convex cones such that $\calK_1\cap\calK_2$ is a subspace. Then $\calK_1-\calK_2$ is closed. \end{lemma}
\begin{IEEEproof} Let $A=\begin{bmatrix} I & -I\end{bmatrix}$ and $\calC= \calK_1 \times \calK_2$. As a consequence of \cite[Thm. 9.1]{Rockafellar:70}, if every $z\in \calC$ such that $Az=0$ belongs to $\lin\calC$, then $A\calC$ is closed. Note that  $z\in \calC$ with $Az=0$ if and only if $z=\begin{bmatrix} y \\ y \end{bmatrix}$ where $y\in \calK_1\cap\calK_2$, proving the lemma.\end{IEEEproof}

We end this subsection with an observation on sequences of nested convex cones.
\begin{lemma}\label{lemm:FGisFD} Let $\calC_\ell\subseteq\R^n$ be a sequence of nested convex cones, i.e. $\calC_\ell\subseteq \calC_{\ell+1}$ such that $\bigcup\limits_{\ell=1}^\infty\calC_\ell=\R^n$. Then, there exists an integer $q$ such that $\calC_q=\R^n$. 
	\end{lemma}
\begin{IEEEproof} Let $\calS$ be a basis of $\R^n$, i.e. $\calS\subseteq\R^n$ is a finite set with $\Lin(\calS)=\R^n$. Then, $-\calS\cup\calS$ is a finite set and hence contained in $\calC_q$ for some integer $q$ since $\bigcup\limits_{\ell=1}^\infty\calC_\ell=\R^n$. Since the only convex cone containing $-\calS\cup\calS$ is $\R^n$, this means that $\calC_q=\R^n$.\end{IEEEproof}

\subsection{Convex processes}
If $H$ is a convex process, then the sets $\dom(H)$, $\im(H)$, and $H(0)$ are convex cones. In addition, we have
\beq\label{eq:conv proc inclusion}
H(x)+H(y)\subseteq H(x+y)
\eeq
for all $x,y\in\dom(H)$. An immediate consequence of this inclusion is that 
\begin{equation} \label{eq: H(x) + H(0)= H(x)}
	H(x)=H(x)+H(0).
\end{equation}
Further, we define the image of a set $\calS$ under $H$ by: 
\begin{equation} \label{eq:image of set}
	H(\calS) = \bigcup_{x\in\calS} H(x).
\end{equation}
This means that for sets $\calS,\calT$, we have:
\bse\label{eq:con pro intersection and sum}
\begin{align}
	H(\calS\cap\calT)&\subseteq H(\calS)\cap H(\calT), \\  H(\calS+\calT) &\subseteq  H(\calS)+H(\calT).
	\end{align}
\ese 	
	
	\subsection{Strong and weak invariance}	
Next, we introduce two notions of invariance under convex processes.

	\begin{dfn}\label{def:invariance}
		Let $H:\mathbb{R}^n\rightrightarrows \mathbb{R}^n$ be a convex process and $\calC\subseteq\R^n$ be a convex cone. We say that $\calC$ is
		\begin{itemize}
			\item \emph{weakly} ${H}$ \emph{invariant} if $H(x)\cap \calC \neq \emptyset$ for all $x\in \calC$.
			\item \emph{strongly} ${H}$ \emph{invariant} if $H(x)\subseteq \calC$ for all $x\in \calC$.
		\end{itemize}
	\end{dfn}
Using the notation from \eqref{eq:image of set}, we can conclude that:
	
	\begin{lemma}\label{lemm:invariance properties.1}
		Let $H:\mathbb{R}^n\rightrightarrows \mathbb{R}^n$ be a convex process. A set $\calW$ is weakly ${H}$ invariant if and only if $\calW\subseteq H^{-1}(\calW)$. A set $\calS$ is strongly ${H}$ invariant if and only if $H(\calS)\subseteq \calS$.
	\end{lemma}
	
	These two invariance notions enjoy the following properties.
	\begin{lemma}\label{lemm:invariance properties} Let $H$ be a convex process. If $\calW$ and $\calS$ are, respectively, weakly and strongly ${H}$ invariant, then $\calW\cap \calS$ and $\calS-\calW$ are, respectively, weakly and strongly ${H}$ invariant.
	\end{lemma}
	
	\begin{IEEEproof} To prove the first part of the statement, let $x\in\calW\cap\calS$. Since $\calW$ is weakly invariant, there exists $y\in H(x)\cap \calW\nonempty$. In view of strong invariance of $\calS$, we have $H(x)\subseteq\calS$. Therefore, $y\in H(x)\cap\calS\cap\calW\nonempty$ and hence $\calW\cap \calS$ is weakly $H$ invariant.
	
	For the second part, let $x\in\calS-\calW$. Then, there exists $s\in\calS$ and $w\in\calW$ such that $x=s-w$. If $x\not\in\dom(H)$, we have $H(x)=\emptyset\subseteq\calS-\calW$. Suppose that $x\in\dom(H)$. Since $\calW$ is weakly invariant, $w\in\dom(H)$. As $\dom(H)$ is a convex cone, we see that $s=x+w\in\dom(H)$. Note that $H(x)+H(w)\subseteq H(s)\subseteq\calS$ since $H$ is a convex process and $\calS$ is strongly invariant. Weak invariance of $\calW$ implies that there exists $z\in H(w)\cap\calW\nonempty$. Then, we have $H(x)+\pset{z}\subseteq \calS$ which implies that $H(x)\subseteq\calS-\pset{z}\subseteq\calS-\calW$. Consequently, $\calS-\calW$ is strongly $H$ invariant.\end{IEEEproof}

\subsection{Feasible, reachable, null-controllable sets}
	If $H$ is a convex process it is immediate that the feasible, reachable and null-controllable set are all convex cones. By inspection, it is immediate that $\calN(H)$ is weakly $H$ invariant, $\calF(H)$ is the largest weakly $H$ invariant cone and $\calR(H)$ is the smallest strongly invariant cone. 
	
	In addition to their definitions in \eqref{eq:reach null def}, it is easy to show that the reachable and null-controllable sets are equal to	
	\[\calR(H) = \bigcup_{\ell=0}^\infty H^{\ell}(0) \qand \calN(H) = \bigcup_{\ell=0}^\infty H^{-\ell}(0). \]	
	This shows that $\calR(H)=\calN(H^{-1})$ and vice-versa. As a consequence, we see that $\calR(H)$ is weakly $H^{-1}$ invariant and $\calN(H)$ strongly $H^{-1}$ invariant.	
	
	\blem\label{lemm:RH invariance}
	Let $H$ be a convex process and let $\calR=\calR(H)$. If $\dom(H)-\calR$ is a subspace, then $\calR-\calR=\calR_+$.
	\elem
	
\begin{IEEEproof} It follows from \eqref{eq:lh-+} that $H^\ell(0)\subseteq L_+^\ell(0)$ for all $\ell\in\N$. Then, we get $\calR-\calR\subseteq \calR_+-\calR_+=\calR_+$ where the last equality follows from the fact that $\calR_+$ is a subspace. As $\calR_+$ is the smallest strongly $L_+$ invariant cone, the reverse inclusion $\calR_+\subseteq\calR-\calR$ would follow if $\calR-\calR$ is strongly $L_+$ invariant. Therefore, it suffices to show that $L_+(\calR-\calR)\subseteq \calR-\calR$. Let $x\in\calR-\calR$. Then, there exist $r_1,r_2\in\calR$ such that $x=r_1-r_2$. If $x\not\in\dom(L_+)=\dom(H)-\dom(H)$, we readily have $L_+(x)=\emptyset\subseteq\calR-\calR$. Suppose that $x\in\dom(H)-\dom(H)$ and $y\in L_+(x)$. From the definition of $L_+$ \eqref{eq:def of L- and L+}, it follows that there exist $x_1,x_2\in\dom(H)$ such that $x=r_1-r_2=x_1-x_2$ and $y\in H(x_1)-H(x_2)$. Note that $r_1-x_1=r_2-x_2\in\calR-\dom(H)$. Since $\dom(H)-\calR$ is a subspace, we have $\calR-\dom(H)=\dom(H)-\calR$. Therefore, there exist $\barx\in\dom(H)$ and $\barr\in\calR$ such that $r_1-x_1=r_2-x_2=\barx-\barr$. Since $H$ is a convex process, we have $H(\barx)+H(x_i)\subseteq H(\barx+x_i)=H(r_i+\barr)$ for $i=1,2$. This leads to $H(x_1)-H(x_2)\subseteq H(r_1+\barr)-H(r_2+\barr)+H(\barx)-H(\barx)$ since $\barx\in\dom(H)$. From $H(\calR)\subseteq\calR$, we know that $H(r_1+\barr)-H(r_2+\barr)\subseteq\calR-\calR$ since $r_1,r_2,\barr\in\calR$. Thus, it suffices to show that $H(\barx)-H(\barx)\in\calR-\calR$ for all $\barx\in\dom(H)$. Let $\barx\in\dom(H)$. As $0\in\calR$, we have $\barx\in\dom(H)-\calR$. Since $\dom(H)-\calR$ is a subspace, $-\barx\in\dom(H)-\calR$ and hence there exist $\xi\in\dom(H)$ and $\eta\in\calR$ such that $-\barx=\xi-\eta$. This means that $H(\barx)+H(\xi)\subseteq H(\barx+\xi)=H(\eta)$ as $H$ is a convex process. Since $\xi\in\dom(H)$, we get $-H(\barx)\subseteq H(\xi)-H(\eta)$. Then, we obtain $H(\barx)-H(\barx)\subseteq H(\barx)+H(\xi)-H(\eta)\subseteq H(\barx+\xi)-H(\eta)=H(\eta)-H(\eta)$. As $H(\eta)\subseteq \calR$, we finally arrive at $H(\barx)-H(\barx)\subseteq \calR-\calR$ since $\eta\in\calR$.\end{IEEEproof}

	\subsection{Dual processes}
	
	In terms of the graph, we can write definition \eqref{eq:def of dual} of the negative dual as:
	\begin{equation} \label{eq:def of dual graph}
	\graph (H^-) = \begin{bmatrix} 0 & I \\ -I & 0 \end{bmatrix} \big( \graph (H)\big)^-
	\end{equation} 
	and similar for the positive dual $H^+$. 
	\begin{lemma}\label{lemm:dualprops} Let $H$ be a convex process, then the following hold: 
		\begin{enumerate}
		    \item\label{lemmitem:dualprops0} $\gr(H^-)=-\gr(H^+)$.
			\item\label{lemmitem:dualprops1} $(H^{-1})^- =(H^+)^{-1} $.
			\item\label{lemmitem:dualprops2} $(\dom H)^- = -H^-(0)=H^+(0)$. 
		\end{enumerate}
	\end{lemma}
	\begin{IEEEproof} Statements \eqref{lemmitem:dualprops0} and \eqref{lemmitem:dualprops1} follow immediately from \eqref{eq:def of dual graph} and the fact that 
	\[\graph (H^{-1}) = \begin{bmatrix} 0 & I \\ I &0 \end{bmatrix} \graph (H).\] 
	Lastly, \eqref{lemmitem:dualprops2} is immediate from $\eqref{eq:def of dual}$ by taking $q=0$.\end{IEEEproof}
	
	The image of a convex cone under a convex process enjoys the following duality relation which is slight generalization of \cite[Thm 2.5.7]{Aubin:90}. 
	
	\begin{proposition}\label{prop:duality of images}
		Let $H : \mathbb{R}^n \rightrightarrows \mathbb{R}^n$ be a convex process and $K$ be a convex cone such that
		$\dom(H)-K$ is a subspace. Then, 
		$$
		\big(H(K)\big)^- =(H^-)^{-1}( K^-).
		$$
	\end{proposition}
	\begin{IEEEproof} We can write the left-hand side in terms of the graph by:
	\begin{align*}
	(H(K))^- &= \left(\begin{bmatrix} 0 & I\end{bmatrix} \left(\gr(H)\cap K\times \mathbb{R}^n\right)\right)^-.\\
	\intertext{Using \cite[Cor. 16.3.2]{Rockafellar:70}, this means that:}
	(H(K))^- &= \begin{bmatrix} 0 \\ I\end{bmatrix}^{-1} \left(\gr(H)\cap K\times \mathbb{R}^n\right)^-\\
	&= \begin{bmatrix} 0 \\ I\end{bmatrix}^{-1} \cl\left(\gr(H)^- + K^- \times \{0\}\right)\\
	&= \begin{bmatrix} 0 \\ I\end{bmatrix}^{-1} \begin{bmatrix} 0 & I\\ -I & 0\end{bmatrix}^{-1}\cl\left(\gr(H^-) + \{0\}\times K^+\right)
	\end{align*} 
	By Lemma~\ref{lemm:dualprops}.\ref{lemmitem:dualprops2} we have
	$$[\dom(H)-K]^+ = H^-(0)\cap K^- $$
	Therefore, from our assumption follows that $\graph (H^-) \cap \{0\}\times K^-$ is a subspace, which allows us to use Lemma~\ref{lemm:sum of closed sets} to reveal that
$\graph (H^-) + \{0\}\times K^+$ is closed. 
	Therefore we can drop the closure from our derivation, and find:
	\begin{align*}
	(H(K))^- &\!=\! \begin{bmatrix} I \\ 0\end{bmatrix}^{-1}\!\!\!\!\!\left(\gr(H^-) + \{0\}\times K^+\right)\!=\! (H^-)^{-1}(K^-).
	\end{align*}
	Thus proving the lemma.\end{IEEEproof}
	
	The previous result allows us to reveal a relation between duality and invariance. 
	
	\begin{theorem}\label{thm:strong inv becomes weak in the dual}
	Let $H : \mathbb{R}^n \rightrightarrows \mathbb{R}^n$ be a convex process and $K$ be a convex cone such that $K-\dom(H)$ is a subspace. Then, $K^-$ is weakly $H^-$ invariant if $K$ is strongly $H$ invariant. Conversely, $K$ is strongly $H$ invariant if $K$ is closed and $K^-$ is weakly $H^-$ invariant. 
	\end{theorem}
	
	\begin{IEEEproof} Suppose that $K$ is strongly $H$ invariant. Then, $H(K) \subseteq K$ in view of Lemma~\ref{lemm:invariance properties.1}. Hence, $K^- \subseteq [H(K)]^-$. From Proposition~\ref{prop:duality of images}, we have $K^- \subseteq (H^-)^{-1}(K^-)$. Therefore, $K^-$ is weakly $H^-$ invariant due to Lemma~\ref{lemm:invariance properties.1}.
	
	Now suppose that $K$ is closed and $K^-$ is weakly $H^-$ invariant. Then, $K^- \subseteq (H^-)^{-1}(K^-)$. Proposition~\ref{prop:duality of images} implies that $\big([H(K)]^-\big)^- \subseteq (K^-)^- $. This means that $\cl(H(K)) \subseteq K$ since $K$ is closed. Hence, $H(K)\subseteq K$. In other words, $K$ is strongly $H$ invariant due to Lemma~\ref{lemm:invariance properties.1}.\end{IEEEproof} 	

	There is a link between weakly $H$ invariant cones and eigenvalues of the dual of $H$, given in \cite[Thm. 3.2]{j37}. 
		
		\begin{proposition}\label{thm:eigthm}
			Let $H:\mathbb{R}^n\rightrightarrows\mathbb{R}^n$ be a closed convex process and $\calK\neq \zset$ be a closed convex cone in $\R^n$. Suppose that $\calK$ does not contain a line, $H(0)\cap \calK=\{0\}$, and $\calK$ is weakly $H$ invariant. Then, $\calK$ contains an eigenvector of $H$ corresponding to a nonnegative eigenvalue.
		\end{proposition}
	
\section{Proofs}\label{sec:proofs}
\subsection{Proof of Theorem~\ref{thm:reach}} 
We begin with an alternative characterization of reachability. 
\begin{lemma}\label{lemm:char of reach} Let $H$ be a convex process such that $\dom H +\rrm =\R^n$. Then, $H$ is reachable if and only if $\calR(H)=\R^n$. \end{lemma}
	\begin{IEEEproof} To prove sufficiency, we assume that $\calF(H)\subseteq\calR(H)$. By \cite[Lemma 4.3]{j37} we know that if $\dom H +\rrm =\R^n$, then $\calF(H)+\rrm = \R^n$. By our assumption, we therefore know that $\R^n=\calF(H)+\rrm \subseteq \calR(H)+\rrm = \calR(H)$. 
	For necessity, suppose that $\calR(H)=\R^n$. Then, clearly $\calF(H)\subseteq\calR(H)$.\end{IEEEproof} 

Next, we relate the eigenvectors of the dual process $H^-$ to the reachable set $\calR(H)$.

\blem\label{l:eig H dual vs R(H)}
Let $H$ be a convex process. If $\xi$ is an eigenvector of $H^-$ corresponding to a nonnegative eigenvalue, then $\xi\in\calR(H)^-$.
\elem

\begin{IEEEproof} Let $\lambda\xi \in H^-(\xi)$ for some $\lambda\geq 0$. Note that $(\lambda^{\ell}\xi,\lambda^{\ell+1}\xi)\in\graph H^-$ for any $\ell\geq 0$. Now take $\eta\in \calR(H)$. Then, $\eta \in H^q(0)$ for some $q$. Hence, there exists a finite sequence $(x_k)_{k=0}^q$ such that $x_0=0$, $x_q=\eta$, and $(x_k,x_{k+1})\in\graph H$ for $k=0,\ldots,q-1$. By the definition of the dual process in \eqref{eq:def of dual}, we know that
$ \langle \lambda^{\ell+1}\xi,x_k \rangle \geq \langle \lambda^\ell\xi, x_{k+1} \rangle  $ 
	for any $\ell\geq 0$ and $k=0,\ldots, q -1$. In particular we can conclude that 
	\[ 0= \langle \lambda^{q}\xi,x_0 \rangle \geq \cdots \geq \langle \lambda\xi, x_{q-1}\rangle \geq \langle \xi,\eta\rangle . \]
This shows that $\xi\in\calR(H)^-$.\end{IEEEproof}
	
Based on these result, we can prove Theorem~\ref{thm:reach} as follows. For the necessity, suppose that $H$ is reachable. Then, it follows from Lemma~\ref{lemm:char of reach} that $\calR(H)=\R^n$. Hence, we see that $\calR_+=\R^n$ since $\calR(H)\subseteq\calR_+$. In addition, we have $\calR(H)^-=\zset$. Therefore, Lemma~\ref{l:eig H dual vs R(H)} implies that $H^-$ has no nonnegative eigenvalues. 

For the sufficiency, suppose that $\rrp=\R^n$ and $H^-$ has no nonnegative eigenvalues. Suppose, on the contrary, that $H$ is not reachable. Then, Lemma~\ref{lemm:char of reach} implies that $\calR(H)\neq\R^n$. Since $\Lin(\calR(H))=\calR_+=\R^n$ due to Lemma~\ref{lemm:RH invariance} and the hypothesis, we see that $\calR(H)$ does not contain a line. From the hypothesis $\dom H+\calR_-=\R^n$ and the fact that $\calR_-\subseteq\calR(H)$, we see that $\dom H-\calR(H)=\R^n$. Then, Theorem~\ref{thm:strong inv becomes weak in the dual} implies that $\calR(H)^-$ is weakly $H^-$ invariant since $\calR(H)$ is strongly $H$ invariant. Using the properties of the dual from \eqref{eq:dual of sum and intersection} and Lemma~\ref{lemm:dualprops}.\ref{lemmitem:dualprops2}, we see that	
	\[ H^-(0)\cap\calR(H)^- = (\dom H -\calR(H))^+ = \pset{0}\]
Now, Proposition~\ref{thm:eigthm} applied to the convex process $H^-$ and the cone $\calR(H)^-$ implies that there exists an eigenvector in $\calR(H)^-$ corresponding to a nonnegative eigenvalue of $H^-$. This is a contradiction, proving that $\calR(H)=\R^n$.

	\subsection{Proof of Theorem~\ref{thm:nullc}} 
Similar to the proof of the reachability theorem, we begin with an alternative characterization of null-controllability.	
	\begin{lemma}\label{l:null-cont alt}
Let $H$ be a convex process such that $\dom H+\rrm =\R^n$. Then, $H$ is null-controllable if and only if $\calN(H)-\calR(H)=\R^n$.  
	\end{lemma}
	\begin{IEEEproof} To prove necessity, suppose that $H$ is null controllable, that is $\calF(H)=\calN(H)$. By \cite[Lemma 4.3]{j37} we know that if $\dom H +\rrm =\R^n$, then $\calF(H)+\rrm = \R^n$. By the hypothesis, we therefore get
	\[ \R^n =\calF(H)+\rrm =\calN(H)+\rrm \subseteq \calN(H)-\calR(H).\] 
Hence, we see that $\calN(H)-\calR(H)=\R^n.$

For sufficiency, suppose that $\calN(H)-\calR(H)=\R^n$. By taking $\calC_\ell = \calN(H)-H^\ell(0)$ and applying Lemma~\ref{lemm:FGisFD}, we see that there exists an integer $q$ such that $\calN(H)-H^q(0)=\R^n$. Now, let $x\in\calF(H)$. Then, there exists $y\in H^q(x)$. Since $\calN(H)-\calR(H)=\R^n$, $y= w-z$ where $w\in \calN(H)$ and $z\in H^q(0)$. From \eqref{eq: H(x) + H(0)= H(x)}, we see that $w \in H^q(x)$. As $\calN(H)$ is strongly $H^{-1}$ invariant, we can conclude that $x\in \calN(H)$. Consequently, we obtain $\calF(H)\subseteq\calN(H)$ and hence $H$ is null-controllable.\end{IEEEproof}

We proceed with the proof of Theorem~\ref{thm:nullc} as follows. For necessity, suppose that $H$ is null-controllable. Therefore, we have $\calN(H)-\calR(H)=\R^n$ due to Lemma~\ref{l:null-cont alt}. Let $\lambda\xi \in H^-(\xi)$ for some $\lambda> 0$ and $\xi\in\R^n$. We already know from Lemma~\ref{l:eig H dual vs R(H)} that
\beq\label{eq:xi in dual of RH}
\xi\in\calR(H)^-.
\eeq
From Lemma~\ref{lemm:dualprops}.\ref{lemmitem:dualprops1}, we know that $(H\inv)^-=(H^+)\inv$. Since $\lambda>0$, we have $\lambda\inv(-\xi)\in(H\inv)^-(-\xi)$ due to Lemma~\ref{lemm:dualprops}.\ref{lemmitem:dualprops0}. Therefore, Lemma~\ref{l:eig H dual vs R(H)} implies that
\beq 
-\xi\in \calR(H\inv)^-=\calN(H)^-.
\eeq 
This, together with \eqref{eq:xi in dual of RH}, results in $\xi\in\calR(H)^-\cap\calN(H)^+$. Note that $\calR(H)^-\cap\calN(H)^+=(\calR(H)-\calN(H))^-=(\R^n)^-=\zset$. Therefore, we can conclude that $H^-$ has no positive eigenvalues.  

For sufficiency, suppose that $H^-$ has no positive eigenvalues but $H$ is not null-controllable. Then, we know from Lemma~\ref{l:null-cont alt} that $\calR(H)-\calN(H)\neq \R^n$. Also, we know from Lemma~\ref{lemm:invariance properties} that $\calR(H)-\calN(H)$ is strongly $H$ invariant since $\calR(H)$ and $\calN(H)$ are, respectively, strongly and weakly invariant. Note that the hypothesis $\dom H-\calR_-=\R^n$ implies that $\dom H-(\calR(H)-\calN(H))=\R^n$ since $\calR_-\subseteq\calR(H)$ and $0\in\calN(H)$. As such, we can apply Theorem~\ref{thm:strong inv becomes weak in the dual} to conclude that $(\calR(H)-\calN(H))^-$ is weakly $H^-$ invariant. From Lemma~\ref{lemm:RH invariance}, we see that $\Lin(\calR(H))=\R^n$. This means that 
$\calR(H)^-$ does not contain a line. Since the cone $(\calR(H)-\calN(H))^-$ is contained in $\calR(H)^-$, it cannot contain a line either. Now, observe that 
\begin{align*} H^-(0)\cap(\calR(H)-\calN(H))^- &\subseteq (\dom H -\calR(H))^+ \\
\subseteq (\dom H +\rrm)^+ &= (\R^n)^+ = \{0\} .
\end{align*}
This means we can apply Proposition~\ref{thm:eigthm} to the convex process $H^-$ and the cone $(\calR(H)-\calN(H))^-$ to conclude that there exists an eigenvector $\xi\in (\calR(H)-\calN(H))^-$ corresponding to a nonnegative eigenvalue $\lambda$ of $H^-$. Since $H^-$ has no positive eigenvalues, $\lambda$ must be zero. Therefore, we have $0 \in H^-(\xi)$, or equivalently $\xi\in (H^-)^{-1}(0)$. By noting that $\im H=\dom H\inv$ and using Lemma~\ref{lemm:dualprops}, we can get $(\im H)^-=(H^-)\inv(0)$. Hence, we see that $\xi\in (\im H)^-$. Recall that $\xi\in(\calR(H)-\calN(H))^-=\calR(H)^-\cap\calN(H)^+$. This means that
\begin{align*}
\xi\in (\im H)^-\cap\calN(H)^+&=(\im H -\calN(H))^-\\
\subseteq (\im H + N_-)^- &=(\R^n)^-=\{0\}.
\end{align*}	
This clearly contradicts with the fact that $\xi$ is an eigenvector. Therefore, we must have $\calR(H)-\calN(H) =\R^n$.

\section{Conclusions}\label{sec:conc}
We have studied reachability and null-controllability for a class of discrete-time systems that are given in the form of difference inclusions with convex processes. Under mild conditions on the domain of a given convex process, we established necessary and sufficient conditions for reachability. For null-controllability, we provided also necessary and sufficient conditions under mild assumptions on both  domains and the images. The results on reachability generalize all existing similar results whereas the results on null-controllability appear for the first time, to the best of our knowledge, in the literature. Moreover, all assumptions we made as well as the conditions we presented can be verified in finite steps.

Future work consists of extending the results from the conic case to the more general case of convex constraints. Also, the framework presented may lead to a similar spectral characterization of the stabilizability problem.

\bibliography{convex_processes}
\bibliographystyle{IEEEtran}

\end{document}